\documentclass[11pt]{amsart}

\usepackage{amsmath,amsthm, amscd, amssymb, amsfonts}

\newcommand{\Hc}{{\mathcal H}}
\newcommand{\otk}{{\otimes_{\ku}}}

\newcommand{\nic}{{\mathfrak B}}

\newcommand{\jc}{{\mathcal J}}

\newcommand{\ot}{{\otimes}}

\newcommand{\Ac}{{\mathcal A}}

\newcommand{\Do}{{\mathcal D}}

\newcommand{\YD}{{\mathcal YD}}

\newcommand{\ku}{{\Bbbk}}

\newcommand{\Z}{{\mathbb Z}}
\newcommand{\Na}{{\mathbb N}}

\renewcommand{\_}[1]{\mbox{$_{\left( #1 \right)}$}}

\newcommand{\id}{\mbox{\rm id\,}}

\theoremstyle{plain}

\numberwithin{equation}{section}

\newtheorem{question}{Question}[section]

\newtheorem{teo}{Theorem}[section]

\newtheorem{lema}[teo]{Lemma}

\newtheorem{cor}[teo]{Corollary}

\newtheorem{prop}[teo]{Proposition}

\theoremstyle{definition}

\newtheorem{defi}[teo]{Definition}

  \newtheorem{exa}[teo]{Example}

\theoremstyle{remark}

\newtheorem{rmk}[teo]{Remark}

\def\pf{\begin{proof}}

\def\epf{\end{proof}}

\theoremstyle{remark}

\begin{document}

\title[Families of  Hopf algebras with the Chevalley property]{Families of finite-dimensional Hopf algebras with the Chevalley property}
\author[Mombelli]{
Mart\'\i n Mombelli }
\thanks{This work was partially  supported by Ministerio de Ciencia
y Tecnolog\'\i a (C\'ordoba),
 Secyt (UNC), CONICET, Argentina}
\address{ FAMAF-CIEM (CONICET),
Universidad Nacional de C\'ordoba,\newline Medina Allende s/n,
(5000) Ciudad Universitaria, C\'ordoba, Argentina}
\email{mombelli@mate.uncor.edu, martin10090@gmail.com
\newline \indent\emph{URL:}\/ http://www.mate.uncor.edu/$\sim$mombelli}
\begin{abstract} We introduce and study new families of finite-dimensional Hopf algebras with the Chevalley property that are not pointed nor semisimple arising as twistings of quantum linear spaces. These Hopf algebras generalize the examples introduced in \cite{AEG}, \cite{EG1}, \cite{EG2}.
\end{abstract}

\subjclass[2000]{16W30, 18D10, 19D23}

\date{\today}
\maketitle

\section{Introduction}

A \emph{Drinfeld's twist}, or simply a \emph{twist}, for a Hopf algebra $H$ is an invertible element $J\in H\ot H$ satisfying a certain non-linear equation, which in some sense is dual to the notion of a 2-cocycle. Any twist gives rise to a new Hopf algebra $H^J$ constructed over the same underlying algebra $H$ such that if $R$ is a quasitriangular structure for $H$ then $R^J=J_{21}RJ^{-1}$ is a quasitriangular structure for $H^J$. The twisting procedure is a very powerful tool to construct new examples of (quasitriangular) Hopf algebras from old and well-known ones.

It has been used in \cite{N}, \cite{GN1} to construct simple Hopf algebras by twisting group algebras, see also \cite{GN2}. In  \cite{AEG} the authors introduce some families of triangular Hopf algebras by twisting finite supergroups. The idea of the construction is roughly the following. They begin with a finite group $G$  with a central element of order two $u\in G$,  $V$  a finite-dimensional representation of $G$ where $u$ acts as $-1$ and $B\in S^2(V)$. The tensor product of the exterior algebra $\wedge(V)$ and the group algebra $\ku G$ is a braided Hopf algebra in the category of Yetter-Drinfeld modules over $\Z_2$. The twisting of $\wedge(V)\ot \ku G$ by the exponential $e^B$ of the symmetric element $B$ gives a new braided Hopf algebra and by bosonization a (usual) Hopf algebra.

This idea was developed further in a series of papers  \cite{EG2}, \cite{EG3}, \cite{EG4} culminating in the classification of finite-dimensional triangular Hopf algebras. Also in \cite{EG1} further properties of this family have been studied.

The main goal of this paper is generalize this construction replacing the group $\Z_2$ by an arbitrary finite
 Abelian group $\Gamma$ and $V$ by a quantum linear space over $G$ such that $\Gamma$
 is in the center of $G$. As expected, the role of the exterior algebra is played by the corresponding
 Nichols algebra $\nic(V)$. This construction gives in some cases new examples of finite-dimensional
Hopf algebras.

The contents of the paper are the following. In section \ref{s1} we briefly recall the definition of quantum linear
 space over a finite group and the definition of twist over a braided Hopf algebra. 
We associate to any twist over a braided Hopf algebra a twist over the corresponding bosonization.
 For any quantum linear space $V$ we construct families of twists in $\nic(V)$ using the quantum
 version of the exponential map.

In section \ref{s2} for any quantum linear space $V$ over a finite group $G$ and a subspace $W\subseteq V$,
 a subgroup $F\subseteq G$ and a twist $J$ for the group algebra $\ku F$  we introduce the definition of
 the Hopf algebras $A(V,G, W, F, J, \Do)$ and we study the particular class of these Hopf algebras when 
$V=W$, $F=\Gamma$, $J=1\ot 1$, that we denote by $A(V,G,\Do)$. In the non-trivial cases these finite-dimensional
Hopf algebras are new.
We describe some isomorphisms of the Hopf algebras $A(V,G,\Do)$ and we study the algebra structure
 of $A(V,G,\Do)^*$ from which we give necessary and sufficient conditions for the Hopf algebra $A(V,G,\Do)$ to be pointed.

\subsection*{Acknowledgments.} This work was  mainly done during a visit to the ICB, Universidad de Cuyo,
 Mendoza, Argentina. I would like to express my gratitude to S. Simondi and Y. Gonzalez for the warm hospitality.
I also want to thank N. Andruskiewitsch for his interest in this work and for his comments that led to a
better presentation of the paper.

\section{Preliminaries and notation}

Throughout the paper $\ku$ will denote an algebraically closed field of characteristic 0 and all vector spaces and algebras are assumed to be over $\ku$.

If $G$ is a group and $(V,\delta)$ is a left $G$-comodule we shall denote $V_g=\{v\in V: \delta(v)=g\ot v\}$ for all $g\in G$.

\medbreak

If $q\in \ku$, $n\in \Na$ denote $(n)_q=1+q+\dots+ q^{n-1}$ and as usual the q-factorial numbers: $n!_q=(1)_q (2)_q\dots (n)_q$. The quantum Gaussian coefficients are defined for $0\leq k \leq n$ by
$${\binom{n}{k}}_{q}=\frac{n!_q}{(n-k)!_q\; k!_q}.$$

The following technical result will be needed later.
\begin{lema} Let $a,i,j,N\in \Na$,  $q\in \ku$ such that $q^N=1$, $1<N$. If  $0\leq a,i,j < N$ and $i+j=N+a$, then
\begin{align}\label{q-binomial-f} \sum_{k=0}^a\; q^{k(k-j)}\, {\binom{j}{k}}_{q}{\binom{i}{a-k}}_{q} =1
\end{align}

\end{lema}

\pf Let $x, y$ be elements in an algebra such that $yx=q\, xy$. Equation \eqref{q-binomial-f} follows by using the quantum binomial formula to expand $(x+y)^a(x+y)^N $ and $(x+y)^i (x+y)^j$, and then compute the corresponding coefficient of the term  $x^a y^N$.
\epf
If $x$ is an element in an algebra such that $x^N=0$ and $q\in\ku$, the \emph{$q$-exponential map } \cite{K} is defined by
$$ \exp_q(x)=\sum_{n=0}^{N-1}\; \frac{1}{(n)!_q}\;\, x^n.$$
The element $ \exp_q(x)$ is invertible, see \cite[Prop. IV.2.6]{K}, and the inverse is given by
$$ \exp_q(x)^{-1}=\sum_{n=0}^{N-1}\; \frac{(-1)^n \, q^{n(n-1)/2}}{(n)!_q}\; x^n. $$
We shall need the following result.
\begin{lema}\cite[Lemma 3.2]{W}\label{q-co} If $q$ is a $N$-th root of 1, $xy=q\, yx$ and $x^i y^{N-i}=0$ for any $i=0\dots N$ then $\exp_q(x+y)=\exp_q(x)\exp_q(y)$.\qed
\end{lema}

\section{Twists in Quantum linear spaces}\label{s1}

We shall describe some families of twists over quantum linear spaces and describe twists for some (usual) Hopf algebras constructed from braided Hopf algebras.

\subsection{Quantum linear spaces}\label{defi:qs}
We shall recall the definition of quantum linear spaces introduced in \cite{AS1}. Let  $G$ be a finite group and $\theta\in \Na$. Let $g_1,\dots,g_\theta\in G$, $\chi_1,\dots,\chi_\theta\in \widehat{G}$. Denote $q_{ij}=\chi_j(g_i)$, $q_i=q_{ii}$, for any $i,j=1,\dots,\theta$. Let $N_i$ be the order of $q_i$ which is assumed to be finite and $N_i>1$. The collection $(g_1,\dots,g_\theta,\chi_1,\dots,\chi_\theta )$ is a \emph{datum for a quantum linear space} if
\begin{align}\label{qsp1}  g_i h&=h g_i, \quad \chi_i\chi_j=\chi_j \chi_i \quad\text{ for all } \; i,j, \text{ and all } h\in G,\\
    \label{qsp2}  q_{ij}q_{ji} &=1 \quad \text{ for all } \; i\neq j.
\end{align}

We shall denote by $\Gamma$ the Abelian group generated by $\{g_i: i=1,\dots,\theta\}$. Let $V$ be the vector space with basis $\{x_1,\dots,x_\theta\}$. With the following maps $V$ is an object in ${}_{G}^{G}\YD$:
$$\delta(x_i)=g_i\ot x_i, \quad h\cdot x_i=\chi_i(h)\, x_i.$$
We shall denote $V=V(g_1,\dots,
g_{\theta},\chi_1, \dots,\chi_{\theta})$. The associated Nichols algebra \cite{AS2} $\nic(V)$ is the algebra generated by elements $\{x_1,\dots,x_\theta\}$ subject to relations
\begin{align}\label{qls3} x_i^{N_i}=0, \quad x_i x_j= q_{ij}\; x_j
x_i \; \text{ if }  i\neq j.
\end{align}

$\nic(V)$ is a braided Hopf algebra in ${}_{G}^{G}\YD$ with coproduct determined by $\Delta(x_i)=x_i\ot 1+ 1\ot x_i$ for all $i=1,\dots,\theta$. The braided Hopf algebra $\nic(V)$ is called a \emph{quantum linear space}.

 Using the quantum binomial formula we get that for any $i, n$
$$\Delta(x_i^n)=\sum_{k=0}^n\; {\binom{n}{k}}_{q_i}\, x_i^k\ot x_i^{n-k}.$$

\medbreak

There is an isomorphism $\nic(V)^*\simeq \nic(V^*)$ of braided Hopf algebras. For any $i=1,\dots, \theta$ and $0\leq r_i < N_i$ define
$X_1^{r_1}\dots X_\theta^{r_\theta}$ the element in $\nic(V)^*$ determined by
$$\langle X_1^{r_1}\dots X_\theta^{r_\theta},  x_1^{s_1}\dots x_\theta^{s_\theta}\rangle=\begin{cases} \prod_{i=1}^\theta\; (r_i)!_{q_i} \quad \text{ if }  \; r_i=s_i, \text{ for all } i=1,\dots, \theta\\
0 \quad \text{ otherwise.} \end{cases}$$

The braided Hopf algebra  $\nic(V)^*$ is generated by elements $X_1,\dots, X_\theta$ subject to relations \eqref{qls3}. The coproduct is determined by $\Delta(X_i)=X_i\ot 1+ 1\ot X_i$ for all $i=1,\dots,\theta$.
\subsection{Twists in braided Hopf algebras}

Let $\Hc$ be a braided Hopf algebra in the category ${}_{G}^{G}\YD$. A \textit{twist} for $\Hc$ is an invertible element $J\in \Hc\ot \Hc$ such that $\delta(J)=1\ot J$ and
\begin{align}\label{twist} (\Delta\ot\id)(J)(J\ot 1)=(\id\ot\Delta)(J)(1\ot J), \;\; (\varepsilon\ot\id)(J)=1=(\id\ot\varepsilon)(J).
\end{align}
Here $\delta: \Hc\ot \Hc\to \ku G\ot \Hc\ot \Hc$ is the coaction of $\Hc\ot \Hc$ in the category ${}_{G}^{G}\YD$. As for usual Hopf algebras there is a new braided Hopf algebra structure on the vector space $\Hc$ with the same algebra structure and coproduct given by $\Delta^J(h)=J^{-1}\Delta(h) J$, for all $h\in \Hc$. This new braided Hopf algebra is denoted by $\Hc^J.$ See \cite{AEG}.

Two twists $J, \widetilde{J}\in \Hc$ are \emph{gauge equivalent} if there exists an invertible element $c\in \Hc$ such that $\epsilon(c)=1$, $\delta(c)=1\ot c$, $g\cdot c=c$ for all $g\in G$ and
$$ \widetilde{J}=\Delta(c) J (c^{-1}\ot c^{-1}).$$
In this case the map $\phi:\Hc^{\widetilde{J}}\to \Hc^J$, $\phi(h)=chc^{-1}$ is an isomorphism of braided Hopf algebras.

\begin{rmk} The product of elements in $\Hc\ot \Hc$ or in $\Hc\ot \Hc\ot \Hc$, as in equation \eqref{twist}, is the product in the tensor product algebra as an object in the corresponding braided tensor category.
\end{rmk}

The following technical Lemma will be of great use later, the proof is straightforward.

\begin{lema}\label{commut-tw} Let $J, J'$ be twists for $\Hc$. Assume that \begin{align}\label{tw1} (1\ot J)(\id\ot \Delta)(J')=(\id\ot \Delta)(J') (1\ot J),\\
\label{tw2} (J\ot 1)(\Delta\ot \id)(J')= (\Delta\ot \id)(J')(J\ot 1).
\end{align}
Then the product $JJ'$ is a twist for $\Hc$.\qed
\end{lema}

For any element $J\in\Hc\ot\Hc$ we shall denote by $\Hc^*_{(J)}$ the vector space $\Hc^*$ with product  given by

\begin{align}\label{twisted-product} \langle X \ast Y, h\rangle=\langle X, h\_1 (h\_2)\_{-1}\cdot J^1\rangle\langle Y, (h\_2)\_{0} J^2\rangle,
\end{align}
for all $X, Y\in \Hc$, $h\in \Hc$.
\begin{lema}\label{tw-alg} The above product in $\Hc^*_{(J)}$ is associative  if and only if $J$ satisfies
$$(\Delta\ot\id)(J)(J\ot 1)=(\id\ot\Delta)(J)(1\ot J).$$
In particular $\Hc^*_{(J)}$ is an algebra with unit $\varepsilon$ in the category ${}_{G}^{G}\YD$ if and only if $J$ is a twist for $\Hc$.\qed
\end{lema}

Let us mention some applications of Lemma \ref{twisted-product}. Denote $W$ the 1-dimensional vector space generated by $x$. The space $W$ is a Yetter-Drinfeld module over $G$ with structure maps given by
$$\delta(x)=g\ot x, \quad f\cdot x=\chi(f)\, x,$$
where $g\in G$, $\chi\in \widehat{G}$, and $q=\chi(g)$. Assume that $q$ has order $N>1$ and that $g^N=1$. Thus $\nic(W)=\ku[x]/(x^N)$. For any $\xi\in \ku$ denote \begin{align}\label{twist-1p} J_\xi=1\ot 1+\sum_{k=1}^{N-1}\; \frac{\xi}{(N-k)!_q k!_q}\;\, x^{N-k}\ot x^k.
\end{align}
\begin{prop} $J_\xi$ is a twist for $\nic(W)$.
\end{prop}
\pf Clearly $J_\xi$ is invertible with inverse given by $J_{-\xi}$, also $(\varepsilon\ot\id)(J_\xi)=1=(\id\ot\varepsilon)(J_\xi).$ Let us prove that $\nic(W)^*_{(J_\xi)}$ with the product \eqref{twisted-product} is associative. The vector space $\nic(W)^*$ has a basis consisting of elements  $\{X^i:i=0,\dots, N-1\},$ where $\langle X^i, x^j \rangle=\delta_{i,j} \;(i!)_q$ for any $i,j\in \{1,\dots, N-1\}$. For any
$0\leq i,j< N$ we have that
\begin{align}\label{prod1} X^i\ast X^j=\begin{cases} X^{i+j} \quad \text{ if } i+ j< N\\
\xi X^a\quad \text{ if } i+j=N+a, 0\leq a<N.
 \end{cases}\end{align}
Hence $X^i\ast \big(X^k \ast X^j\big)= \big(X^i \ast X^k\big)\ast X^j$ for any $i,j, k$ and the product $\ast$ is associative. Let us prove that $X^i\ast X^j=\xi\, X^a $ when $i+j=N+a$, the other case is straightforward. By definition if $h\in\nic(W)$ then $\langle X^i\ast X^j, h\rangle$ is equal to
$$ \langle  X^i, h\_1 \rangle\langle X^j, h\_2 \rangle+\!\sum_{k=1}^{N-1}\;  \frac{\xi}{(N-k)!_q k!_q}\; \langle X^i, h\_1 h\_{2}\_{-1}\cdot x^{N-k}\rangle \langle X^j,h\_{2}\_{0} x^k\rangle.$$
Thus it is clear that if $h\neq x^a$ then $\langle X^i\ast X^j, h\rangle=0$. We have also that $\langle  X^i, (x^a)\_1 \rangle\langle X^j, (x^a)\_2 \rangle=0$, then $\langle X^i\ast X^j, x^a\rangle$ is equal to
\begin{align*} &=\sum_{k=1}^{N-1}\sum_{l=0}^{a}\;{\binom{a}{l}}_{q}  \frac{\xi}{(N-k)!_q k!_q}\; \langle X^i, x^{a-l} g^l\cdot x^{k}\rangle \langle X^j,x^l x^{N-k}\rangle\\
&=\sum_{k=1}^{N-1}\sum_{l=0}^{a}\;{\binom{a}{l}}_{q}  \frac{\xi\, q^{lk}}{(N-k)!_q k!_q}\; \langle X^i, x^{a-l+k} \rangle \langle X^j, x^{N-k+l}\rangle\\
&=\sum_{l=0}^{a}\;{\binom{a}{l}}_{q}  \frac{\xi\,\, q^{l(i+l-a)}\,\, i!_q\,\, j!_q}{(N-i-l+a)!_q (i+l-a)!_q}\\
&=\xi\, a!_q\, \sum_{l=0}^{a}\; q^{l(l-j)} \; {\binom{j}{l}}_{q} {\binom{i}{a-l}}_{q}= \xi\, a!_q.
\end{align*}
The last equality follows from \eqref{q-binomial-f}. \epf

Let $V=V(g_1,g_2)$ be the 2-dimensional Yetter-Drinfeld module for some  datum of a quantum linear space. Assume that $g_1g_2=1$ and that $N_1=N=N_2$, $q_1=q=q^{-1}_2$ is a $N$-th primitive root of unity. Then $\nic(V)$ is  the algebra generated by $x, y$ subject to relations
$$ x^{N}=0=y^{N}, \quad xy= q\, yx.$$
For any $a\in \ku$ set $B=a\; x\ot y$ and $J_a=\exp_q(B)$. It follows by a straightforward computation that $(1\ot B)(B\ot 1)= (B\ot 1)(1\ot B),$ hence the exponentials commute $\exp_q(1\ot B)\exp_q(B\ot 1)=\exp_q(B\ot 1)\exp_q(1\ot B)$.

\begin{prop} $J_a$ is a twist for $\nic(V)$ in the category ${}_{\Gamma}^{\Gamma}\YD$.
\end{prop}
\pf The proof goes in a similar way to the proof of \cite[Thm. 3.3]{W}. It is immediate to verify that $(\varepsilon\ot\id)(J_a)=1=(\id\ot\varepsilon)(J_a).$ First note that
$$(\Delta \ot \id)(\exp_q(B))=\exp_q((\Delta \ot \id)(B)),\; \;  \exp_q(B)\ot 1= \exp_q(B\ot 1)$$ and  $$( \id\ot \Delta)(\exp_q(B))=\exp_q((\id\ot \Delta)(B)),\; \;  1\ot\exp_q(B) = \exp_q(1\ot B).$$

Let us denote $C=a (x\ot 1\ot y)$, then
$$(\Delta \ot \id)(B)= C + 1\ot B, \quad (\id\ot \Delta)(B)=C+ B\ot 1.$$
It follows easily that $q\;(B\ot 1)C= C (B\ot 1)$, $(1\ot B)C=q\; C (1\ot B)$ and for any $k=0\dots N$ we have that $C^k (B\ot 1)^{N-k}=C^k (1\ot B)^{N-k}=0$. Then using Lemma \ref{q-co} we get that
\begin{align*} (\Delta\ot\id)(J_a)(J_a\ot 1)&=\exp_q((\Delta \ot \id)(B))
\exp_q(B\ot 1)\\
&=\exp_q(C + (1\ot B))
\exp_q(B\ot 1)\\
&=\exp_q(C)\exp_q(1\ot B)\exp_q(B\ot 1)\\
&=\exp_q(C)\exp_q(B\ot 1)\exp_q(1\ot B)\\
&=\exp_q(C+ B\ot 1)\exp_q(1\ot B)\\
&=\exp_q((\id\ot \Delta)(B))\exp_q(1\ot B)=(\id\ot\Delta)(J_a)(1\ot J_a).
\end{align*}\epf

\subsection{A Hopf algebra associated to a  braided Hopf algebra} Let $\Gamma$ be a finite Abelian group and $\Hc\in {}_{\Gamma}^{\Gamma}\YD$ such that $\Gamma$ is a subgroup of the group-like elements in $\Hc$ and for any $g\in \Gamma$, $\delta(g)=1\ot g$. Here $\delta:\Hc\to \ku\Gamma\ot \Hc$ is the coaction. Inspired by \cite{AEG} we shall construct a Hopf algebra $H$ such that the tensor categories of representations of $H$ and $\Hc$ are equivalent.

\medbreak

Consider the bosonization $\Hc\#\ku \Gamma$. The ideal $I$ generated by elements $ h \# 1-  1\# h$ for all $h\in \Gamma$ is a Hopf ideal. Define $H=\Hc\#\ku \Gamma/I$. The class of an element $x\ot g\in \Hc\#\ku \Gamma$ in the quotient $H$ will be denoted by $\overline{x\ot g}$.

\medbreak

If $\jc\in \Hc\ot \Hc$ is a twist, define
\begin{align}\label{from-braided-twist} J=\overline{(\id\ot\delta)(\jc)\# 1}=\overline{\jc^1\#\jc^2\_{-1}\ot \jc^2\_0\#1}.
\end{align}
We shall use the notation $\jc=\jc^1\ot \jc^2=j^1\ot j^2$.
\begin{teo}\label{twists-b} The element $J\in H\ot H$ is a twist. If $\jc$ is $\Gamma$-invariant, that is $g\cdot\jc=\jc$ for all $g\in \Gamma$, then $\Hc^{\jc}\#\ku \Gamma/I\simeq H^J$.
\end{teo}
\pf Clearly $J$ is invertible with inverse $J^{-1}= \overline{(\id\ot\delta)(\jc^{-1})\# 1}$.
Applying $(\id\ot\delta\ot \delta)$ to $(\Delta\ot\id)(\jc)(\jc\ot 1)$ we obtain
\begin{align*}\jc^1\_1 \big(\jc^1\_2\_{-1}\jc^2\_{-2}\big)\cdot j^1&\ot \jc^1\_2\_{-1} j^2\_{-1}\ot\\ &\ot \jc^1\_2\_0\jc^2\_{-1}\cdot j^2\_0\ot \jc^2\_{-1}\ot \jc^2\_0.
\end{align*}
Applying to this element $(\id\ot\id\ot\id\ot\Delta\ot\id)$, multiplying the fourth and second tensorands and using the cocommutativity of $\ku\Gamma$ with obtain
\begin{equation}\label{twq1}\begin{split}\jc^1\_1 \big(\jc^1\_2\_{-2}\jc^2\_{-4}\big)\cdot j^1&\ot \jc^1\_2\_{-1} \jc^2\_{-3}j^2\_{-1}\ot\\ &\ot \jc^1\_2\_0\jc^2\_{-2}\cdot j^2\_0\ot \jc^2\_{-1}\ot \jc^2\_0.
\end{split}
\end{equation}

On the other hand, applying $(\id\ot\delta\ot \delta)$ to $(\id\ot\Delta)(\jc)(1\ot\jc)$ we obtain
\begin{align*} \jc^1\ot \jc^2\_{1}\_{-1}j^1\_{-1}\ot \jc^2\_1\_0\jc^2\_2\_{-2}\cdot j^1\_0\ot \jc^2\_2\_{-1}j^2\_{-1}\ot \jc^2\_2\_0j^2\_0.
\end{align*}
Again, applying to this element $(\id\ot\id\ot\id\ot\Delta\ot\id)$ and multiplying the fourth and second tensorands we obtain
\begin{align*} \jc^1\ot \jc^2\_{1}\_{-1}\jc^2\_{2}\_{-2}j^1\_{-1}j^2\_{-1}\_1&\ot \jc^2\_1\_0\jc^2\_2\_{-3}\cdot j^1\_0\ot\\
&\ot  \jc^2\_2\_{-1}j^2\_{-1}\_2\ot \jc^2\_2\_0j^2\_0.
\end{align*}
Using the cocommutativity of  $\ku\Gamma$ and that $\delta(\jc)=1\ot \jc$ then the above element is equal to
\begin{align*} \jc^1\ot \jc^2\_{1}\_{-1}\jc^2\_{2}\_{-3}&\ot \jc^2\_1\_0\jc^2\_2\_{-2}\cdot j^1\_0\ot\\
&\ot  \jc^2\_2\_{-1}j^2\_{-1}\ot \jc^2\_2\_0j^2\_0,
\end{align*}
and since the coproduct is a morphism of $\ku\Gamma$-comodules the above element is equal to
\begin{align}\label{twq2} \jc^1\ot \jc^2\_{1}&\ot \jc^2\_0\_1\jc^2\_0\_2\_{-2}\cdot j^1\_0\ot  \jc^2\_0\_2\_{-1}j^2\_{-1}\ot \jc^2\_2\_0j^2\_0,
\end{align}

The element $(\Delta\ot\id)(J)(J\ot 1)$ is equal to
\begin{align*}(\Delta\ot\id)&(\overline{\jc^1\#\jc^2\_{-1}\ot \jc^2\_0\#1})(J\ot 1)=\\
&=\overline{\jc^1\_1 \big(\jc^1\_2\_{-2}\jc^2\_{-4}\big)\cdot j^1\# \jc^1\_2\_{-1}\jc^2\_{-3} j^2\_{-1}}\ot \\
&\quad\quad\quad\quad\quad\ot\overline{\jc^2\_2\_0\jc^2\_{-2}\cdot j^2\_0\# \jc^2\_{-1}}\ot\overline{\jc^2\_0\#1},
\end{align*}
and $(\id\ot\Delta)(J)(1\ot J)$ equals
\begin{align*}(\id\ot\Delta)&(\overline{\jc^1\#\jc^2\_{-1}\ot \jc^2\_0\#1})(1\ot J)=\\
&=\overline{\jc^1\# \jc^2\_{1}}\ot \overline{\jc^2\_0\_1\jc^2\_0\_2\_{-2}\cdot j^1\_0\#   \jc^2\_0\_2\_{-1}j^2\_{-1}}\ot\\
&\quad\quad\ot \overline{\jc^2\_2\_0j^2\_0\# 1}.
\end{align*}
Since $\jc$ is a twist then equations \eqref{twq1}, \eqref{twq2} are equal, hence we conclude that $(\Delta\ot\id)(J)(J\ot 1)=(\id\ot\Delta)(J)(1\ot J)$.
 It follows easily that the coproduct of $\Hc^{\jc}\#\ku \Gamma/I$ coincides with the coproduct of  $H^J$.\epf

\subsection{The exponential map in quantum linear spaces}

Let  $G$ be a finite group and $\theta\in \Na$,   $(g_1,\dots,g_\theta,\chi_1,\dots,\chi_\theta )$ be a datum for a quantum linear space and $V=V(g_1,\dots,g_\theta,\chi_1,\dots,\chi_\theta )$.

Let $\Do=\{a_{ij}\in \ku: 1\leq i, j\leq \theta, i\neq j\}\cup\{\xi_i\in \ku:1\leq i\leq \theta\}$ be a family of $\theta^2$ scalars. We shall say that $\Do$ is \emph{compatible} with the quantum linear space $V$ if
\begin{align}\label{even-e} a_{ij}=0\; \text{ if }\; g_ig_j\neq 1, \quad \xi_i=0 \;\text{ if }\; g_i^{N_i}\neq 1.
\end{align}

Let $F\subseteq G$ be a subset. We shall say that the family of scalars $\Do$ is \textit{$F$-invariant} if
\begin{align}\label{conditions1}\chi_i(g) \chi_j(g)\, a_{ij}&= a_{ij}, \\
\label{conditions2} \chi^{N_i}_i(g)\, \xi_i&= \xi_i \quad \text{ for all }  g\in F.
\end{align}

In particular if $\Do$ is compatible with $V$ then it is $\Gamma$-invariant, which amounts to
\begin{align}\label{bilinear:conditions1}  a_{ij}&=0\quad  \text{ if } \; q_{ik}q_{jk}\neq 1 \text{ for some } k=1,\dots,\theta,\\
\label{bilinear:conditions2} \xi_i&=0 \quad  \text{ if } \;\; q^{N_i}_{ij} \neq 1 \text{ for some } j=1,\dots,\theta.
\end{align}
If $g\in G$ and $\Do_1=\{a_{ij}\in \ku: 1\leq i, j\leq \theta, i\neq j\}\cup\{\xi_i\in \ku:1\leq i\leq \theta\}$, $\Do_2=\{a'_{ij}\in \ku: 1\leq i, j\leq \theta, i\neq j\}\cup\{\xi'_i\in \ku:1\leq i\leq \theta\}$ are two families of scalars we will denote
$$\Do_1+\Do_2= \{a_{ij}+a'_{ij}\in \ku\}\cup\{\xi_i+\xi'_i\in \ku:1\leq i\leq \theta\},$$
$$g\cdot \Do_1=\{\chi_i(g) \chi_j(g)\, a_{ij}\in \ku\}\cup\{\chi^{N_i}_i(g)\,\xi_i\in \ku:1\leq i\leq \theta\}.$$
We shall say that a family of scalars $\Do$ is \emph{q-symmetric} if $a_{ij}=-q_{ij} a_{ji}$ for any $ 1\leq i, j\leq \theta, i\neq j$. We shall denote $$\widehat{\Do}=\{b_{ij}\in \ku: 1\leq i< j\leq \theta, i\neq j\}\cup\{\xi_i\in \ku:1\leq i\leq \theta\}$$ where $b_{ij}= q_{ij} a_{ji}-a_{ij}$. Clearly $\widehat{\Do}$ is q-symmetric.

\medbreak

Define $B_{ij}=a_{ij}\; x_i\ot x_j$, $ J_{\xi_i}=1\ot 1+\sum_{k=1}^{N_i-1}\; \frac{\xi_i}{(N_i-k)!_{q_i} k!_{q_i}}\;\, x_i^{N_i-k}\ot x_i^k$ and
\begin{align}\label{definition-twist} J_\Do=\prod^{\theta}_{i=1} J_{\xi_i}\prod_{1\leq i, j\leq \theta, i\neq j} \exp_{q_{ij}}(B_{ij}). \end{align}

\begin{teo} Let $\Do$ be a  compatible family of scalars with $V$. Then $J_\Do$ is a twist for $\nic(V)$  in the category  ${}_{\Gamma}^{\Gamma}\YD$.
\end{teo}
\pf Follows by \eqref{bilinear:conditions1}  that
$$ (1\ot B_{ij})(\id\ot \Delta(B_{kl}))=(\id\ot \Delta(B_{kl}))(1\ot B_{ij}),$$
$$ (B_{ij}\ot 1)(\Delta(B_{kl})\ot\id)=(\Delta(B_{kl})\ot\id) (B_{ij}\ot 1),$$
thus
$$(1\ot \exp_{q_{ij}}(B_{ij}))(\id\ot \Delta(\exp_{q_{kl}}(B_{kl})))= (\id\ot \exp_{q_{kl}}(B_{kl}))(1\ot \exp_{q_{ij}}(B_{ij})),$$
$$ ( \exp_{q_{ij}}(B_{ij})\ot 1)(\Delta(\exp_{q_{kl}}(B_{kl}))\ot\id)=(\Delta(\exp_{q_{kl}}(B_{kl}))\ot\id) ( \exp_{q_{ij}}(B_{ij})\ot 1).$$
Using Lemma \ref{commut-tw} we obtain that $\prod_{1\leq i< j\leq \theta} \exp_{q_{ij}}(B_{ij})$ is a twist. It follows from  \eqref{bilinear:conditions2} that
$$(1\ot x_i^{N_i-k}\ot x_i^k)(\id\ot \Delta)(x_j^{N_j-a}\ot x_j^a)=(\id\ot \Delta)(x_j^{N_j-a}\ot x_j^a)(1\ot x_i^{N_i-k}\ot x_i^k), $$
thus by Lemma \ref{commut-tw}  we conclude that $\prod^{\theta}_{i=1} J_{\xi_i}$ is a twist.

\medbreak

It follows from \eqref{bilinear:conditions2} that $(x_i^{N_i-k}\ot x_i^k)(x_l\ot x_j)=(x_l\ot x_j)(x_i^{N_i-k}\ot x_i^k)$, thus $\exp_{q_{lj}}(B_{lj})$ and $ J_{\xi_i}$ commute. Let $i,j,k=1,\dots,\theta$, $1 \leq a \leq N_k$, then
$(1\ot x_k^{N_k-a}\ot  x_k^{a})(x_i\ot 1 \ot x_j+ x_i\ot x_j\ot 1)$ equals to
\begin{align*} &= x_i\ot x_k^{N_k-a}\ot  x_k^{a} x_j +q_{kj}^a\;  x_i\ot x_k^{N_k-a} x_j\ot  x_k^{a} \\
&=q_{kj}^a\; x_i\ot x_k^{N_k-a}\ot  x_j x_k^{a} + q_{kj}^{N_k}\; x_i\ot x_j x_k^{N_k-a}\ot  x_k^{a}.
\end{align*}
On the other hand $(x_i\ot 1 \ot x_j+ x_i\ot x_j\ot 1)(1\ot x_k^{N_k-a}\ot  x_k^{a})$ equals to
$$ q_{jk}^{N_k-a}\; x_i\ot x_k^{N_k-a}\ot x_jx_k^{a} + x_i\ot x_j x_k^{N_k-a}\ot  x_k^{a}.$$
Using \eqref{bilinear:conditions2} follows that $(1\ot x_k^{N_k-a}\ot  x_k^{a})$ and $(\id\ot\Delta)(x_i\ot x_j)$ commute, hence $(1\ot J_{\xi_i})$ and $(\id\ot\Delta)(\exp_{q_{lj}}(B_{lj}))$ commute. Similarly we can prove that $( J_{\xi_i}\ot 1)$ and $(\Delta\ot\id)(\exp_{q_{lj}}(B_{lj}))$ commute. Using again Lemma \ref{commut-tw} follows that $J_\Do$ is a twist.\epf

\begin{rmk}\label{commutation-t} For  any $1\leq k,l,s,t \leq \theta$ we have that
$$\exp_{q_{kl}}(B_{kl})\exp_{q_{st}}(B_{st})= \exp_{q_{st}}(B_{st})\exp_{q_{kl}}(B_{kl}),$$
$$ \exp_{q_{kl}}(B_{kl})J_{\xi_j}=J_{\xi_j}\exp_{q_{kl}}(B_{kl}),\; J_{\xi_j}J_{\xi_t}=J_{\xi_t}J_{\xi_j}.$$
\end{rmk}

\begin{rmk} It would be interesting to study the exponential map for other types of Nichols algebras.
\end{rmk}

\section{Hopf algebras $A(V,G, W, F, J, \Do)$}\label{s2}

Let  $G$ be a finite group and $\theta\in \Na$,   $(g_1,\dots,g_\theta,\chi_1,\dots,\chi_\theta )$ be a datum for a quantum linear space and $V=V(g_1,\dots,g_\theta,\chi_1,\dots,\chi_\theta )$. As before $\Gamma$ is the Abelian group generated by $\{g_i: i=1,\dots,\theta\}$. Note that $\Gamma$ is contained in the center of $G$. Using ideas contained in \cite{AEG} we shall construct Hopf algebras coming from twisting $\nic(V)\#\ku G$.
\medbreak

By restriction $V$ is an object in ${}_{\Gamma}^{\Gamma}\YD$. The vector space $\nic(V)\otk \ku G$ has a structure of braided Hopf algebra in ${}_{\Gamma}^{\Gamma}\YD$ as follows. The coaction $\delta:\nic(V)\otk \ku G\to \ku\Gamma\otk \nic(V)\otk \ku G$ and the action $\cdot: \Gamma\otk \nic(V)\otk \ku G\to\nic(V)\otk \ku G$ are determined by
$$\delta(v\ot g)=v\_{-1}\ot v\_0\ot g,\quad h\cdot(v\ot g)= h\cdot v\ot g, $$
 for all $v\in \nic(V)$, $g\in G$, $h\in \Gamma$. The product and coproduct in $\nic(V)\otk \ku G$ are given by
$$(v\ot g)(v'\ot g') =v g\cdot v'\ot gg', \quad \Delta(v\ot g)= v\_1\ot g\ot v\_2\ot g,$$
 for all $v,v'\in \nic(V)$, $g,g'\in G$.

 \medbreak

 Let $F$ be a subgroup of $G$ such that $\Gamma\leq F$, let  $W\subseteq V$ be a subspace stable under the action of $F$ and $W_g\subseteq V_g$ for all $g\in F$. In this case we can consider the braided Hopf algebra $\nic(W)\ot \ku F$. Let $\Do$ be an $F$-invariant family of scalars for the quantum linear space $W$ and let $J$ be a twist of $\ku F$.

\begin{lema}  The element $J^1_\Do\ot J^1\ot J^2_\Do\ot J^2$ is a twist for the  braided Hopf algebra $\nic(W)\ot \ku F$ in the category ${}_{\Gamma}^{\Gamma}\YD$.
\end{lema}
\pf Both elements  $J^1_\Do\ot 1\ot J^2_\Do\ot 1$, $1\ot J^1\ot 1\ot J^2$ are twists. Note that for any $g\in F$ we have that $g\cdot J_\Do=J_{g\cdot \Do}$, hence
$$ (1\ot g\ot 1\ot g)(J^1_\Do\ot 1\ot J^2_\Do\ot 1)= J^1_{g\cdot \Do}\ot g\ot J^2_{g\cdot \Do}\ot g.$$ Since $\Do$ is  $F$-invariant follows that equations \eqref{tw1} and \eqref{tw2} are satisfied and by Lemma \ref{commut-tw} $J^1_\Do\ot J^1\ot J^2_\Do\ot J^2$ is a twist.\epf

Abusing  the notation we shall denote by $J_\Do J$ the twist $J^1_\Do\ot J^1\ot J^2_\Do\ot J^2$.

\begin{defi}  Under the above assumptions we define the braided Hopf algebra
$$\Ac(V,G, W, F, J, \Do)=(\nic(V)\ot \ku G)^{J_\Do J}$$ in the category ${}_{\Gamma}^{\Gamma}\YD$. Using the bosonization procedure we construct the Hopf algebra $\Ac(V,G, W, F, J, \Do)\#\ku\Gamma$. The bilateral ideal $I$ generated by elements $1\ot h \# 1- 1\ot 1\# h$ for all $h\in \Gamma$ is a Hopf ideal. Thus we define the Hopf algebra $A(V,G, W, F, J, \Do)$ as the quotient $\Ac(V,G, W, F, J, \Do)\#\ku\Gamma/I$.

If $V=W$, $F=G$ and $J=1\ot 1$ we shall denote the Hopf algebra $A(V,G, V, \Gamma, J, \Do)$ simply by $A(V,G, \Do)$.
\end{defi}

\begin{rmk} Note that if $\Do=0$, that is if $a_{ij}=0=\xi_i$ for all $1\leq i,j\leq\theta$, then $J_\Do=1\ot 1$ and $A(V,G,0)=\nic(V)\#\ku G$.
\end{rmk}

\begin{cor}  $A(V,G, W, F, J, \Do)$ is twist equivalent to $\nic(V)\#\ku G$.
\end{cor}
\pf It follows from Theorem \ref{twists-b}.\epf

\begin{defi}\cite{AEG}  A tensor category is said to have the \emph{Chevalley property} if the tensor product of simple objects is semisimple. A Hopf algebra $H$ has the Chevalley property if the category of left $H$-modules does.
\end{defi}

If $H$ is a Hopf algebra with the Chevalley property and $J\in H\ot H$ is a twist then $H^J$ has the Chevalley property. Hence the families of Hopf algebras $A(V,G, W, F, J, \Do)$ have the Chevalley property.

\begin{exa}\label{exa1}\cite{AEG} Let $G$ be a finite group and $u\in G$ be a central element of order 2. Let $V$ be a $G$-module such that $u\cdot v=-v$ for all $v\in V$. The space $V$ is a Yetter-Drinfeld module over $G$ by declaring $V=V_u$, thus $\Gamma=\Z_2$. Let $\{x_1,\dots, x_\theta\}$ be a basis of $V$. In this case $q_{ij}=-1=q_{ii}$ for all $1\leq i,j\leq\theta$ and the Nichols algebra $\nic(V)$ is the exterior algebra $\wedge V$.

Let $\Do=\{a_{ij}\in \ku: 1\leq i, j\leq \theta, i\neq j\}\cup\{\xi_i\in \ku:1\leq i\leq \theta\}$ be a family of scalars. Note that $\Do$ is automatically $\Z_2$-invariant. Define
$$ B= \sum_{i\neq j} \; a_{ij} \,x_i\ot x_j + \sum_{i=1}^\theta\; \xi_i\, x_i\ot x_i\in V\ot V.$$
Since $q_{ij}=-1$ then $J_\Do= e^B$. Our definition of $A(V,G,\Do)$ coincides with the definition given in \cite{AEG}, see also \cite{EG1}, where this algebra is denoted by $A(V,G,B)$. Note, however, that in \emph{loc. cit.} the authors assume that the element $B$ is symmetric, that is $B\in S^2(V)$.

Let us develop a more particular example.
 Assume that $V$ has a basis $\{x,y\}$. If $a\in\ku$ denote $B_a=a\; x\ot y- a\; y\ot x$.
 In this case the twist $e^{B_a}$ is gauge equivalent to the trivial twist $1\ot 1$. 
Indeed if $c=e^{a\, xy}$ then $e^{B_a}=\Delta(c)(c^{-1}\ot c^{-1}).$
 Thus $A(V,G,B_a)\simeq \wedge V\#\ku G.$ The isomorphism is given by conjugation by $c$, thus
 one can not expect to apply \cite[Prop. 2.1]{EG1} in the general situation.
 Also $\Ac(V,G,B_a)$ is super cocommutative but $B_a$ is not $G$-invariant hence \cite[Corollary 5.3]{EG1} is no longer true when $B$ is not symmetric.
\end{exa}

\begin{exa}\label{1-dim} Let $G$ be a finite group with a character $\chi:G\to \ku^{\times}$. Let be $n\in \Na$ and $\xi\in\ku$. Let $u\in G$ be a central element of order $n$ and $\chi(u)=q$ be a $n$-th primitive root of unity. Let $V$ be the one-dimensional vector space generated by $x$ with structure of Yetter-Drinfeld module over $G$ given by
$$\delta(x)=u\ot x,\quad g\cdot x=\chi(g)\, x,\; \text{ for all } g\in G. $$
The Nichols algebra of $V$ is isomorphic to $\ku[x]/(x^n)$. If $\Do_\xi=\{\xi\}$ then $A(V,G,\Do_\xi)$ is isomorphic to the algebra generated by elements $\{x, g: g\in G\}$ subject to relations
$$ x^n=0,\quad g x= \chi(g)\, xg, \text{ for all } g\in G.$$
The twist in this case is $J_\xi=1\ot 1+ \sum_{k=1}^{n-1}\; \frac{\xi}{(n-k)!_{q} k!_{q}}\;\, x^{n-k}\ot x^k$ and the coproduct is given by formulas
$$\Delta(x)=x\ot 1+ u\ot x, \; \Delta(g)=g\ot g+\sum_{k=1}^{n-1}\; \frac{\xi (\chi^n(g)-1)}{(n-k)!_{q} k!_{q}}\;\, x^{n-k}u^kg\ot x^kg. $$
We shall prove later that if $\xi\neq 0$ and $\chi^n\neq 1$ then $A(V,G,\Do_\xi)$ is not a pointed Hopf algebra and (unfortunately) $A(V,G,\Do_\xi)\simeq A(V,G,\Do_1)$.
\end{exa}

\subsection{Some isomorphisms of $A(V,G,\Do)$} In this section we shall present some isomorphisms of Hopf algebras $A(V,G,\Do)$ in a similar way as in \cite[Prop. 2.1]{EG1}.

\begin{rmk}\label{rmk-coprod} The coproduct of the braided Hopf algebra $\nic(V)\ot \ku G$ is given by:
\begin{align*}\Delta^{J_\Do}(v\ot g)&=J^{-1}_\Do \Delta(v\ot 1)(1\ot g\ot 1\ot g)J_\Do=J^{-1}_\Do J_{g\cdot\Do }  \Delta(v\ot g),
\end{align*}
for all $g\in G, v\in\nic(V)$. This equation follows since $J_\Do \Delta(v\ot 1)=\Delta(v\ot 1)J_\Do$
 for all $v\in\nic(V)$. This implies that if $\Do$ is $G$-invariant then $A(V,G,\Do)\simeq \nic(V)\#\ku G.$ 
In particular if $G=\Gamma$ then $A(V,G,\Do)\simeq \nic(V)\#\ku G.$ 
\end{rmk}

Let $G, G'$ be finite groups, $V, V'$ be quantum linear spaces over $\Gamma$ and $\Gamma'$ respectively with basis $\{x_1,\dots,x_\theta\}$ and $\{x'_1,\dots,x'_\theta\}$ respectively. Let $\Do, \Do'$ be families of scalars such that $\Do$ is $\Gamma$-invariant and $\Do'$ is $\Gamma'$-invariant.

\begin{prop}\label{iso-classes} The Hopf algebras $A(V,G,\Do)$, $A(V',G',\Do')$ are isomorphic provided there is a group isomorphism $\phi:G\to G'$ such that $\phi(\Gamma)=\Gamma'$ and a isomorhism $\eta:V\to \phi^*(V')$ of Yetter-Drinfeld modules over $G$ such that $(\eta\ot \eta)(J_\Do)= J_{\Do'}J_{\widetilde{\Do}}$ where $\widetilde{\Do}$ is $G$-invariant.
\end{prop}
\pf We shall prove that the braided Hopf algebras $ \nic(V)\ot \ku G$, $ \nic(V')\ot \ku G'$ are isomorphic. The map $\eta$ can be extended to an algebra map $\eta:\nic(V)\to \nic(V')$. Define $\psi: \nic(V)\ot \ku G\to \nic(V')\ot \ku G'$ by
$$ \psi(v\ot g)=\eta(v)\ot \phi(g), \text{ for all } v\in \nic(V), g\in G.$$
Clearly $\psi$ is a bijective algebra morphism. If $v\in \nic(V), g\in G$ then
\begin{align*} \Delta^{J_{\Do'}}(\psi(v\ot g))=J^{-1}_{\Do'}J_{g\cdot \Do'} \Delta(\psi(v\ot g)).
\end{align*}
On the other hand
\begin{align*} (\psi\ot \psi)\Delta^{J_{\Do}}(v\ot g)&= (\eta\ot \eta)(J^{-1}_{\Do}J_{g\cdot \Do}) (\psi\ot \psi)\Delta(v\ot g)\\
&= J^{-1}_{\Do'}J^{-1}_{\widetilde{\Do}} J_{\phi(g)\cdot\Do'}J_{\phi(g)\cdot\widetilde{\Do}}(\psi\ot \psi)\Delta(v\ot g)\\
&=  J^{-1}_{\Do'} J_{\phi(g)\cdot\Do'}J^{-1}_{\widetilde{\Do}}J_{\widetilde{\Do}}(\psi\ot \psi)\Delta(v\ot g)\\
&=J^{-1}_{\Do'}J_{g\cdot \Do'} \Delta(\psi(v\ot g)).
\end{align*}
The third equation follows from remark \ref{rmk-coprod}.\epf

\subsection{The algebra structure of $A(V,G,\Do)^*$} In this section we shall describe the algebra structure of $A(V,G,\Do)^*$ following very closely the proof given in \cite{EG1}. As a consequence we shall give necessary and sufficient conditions for the Hopf algebra $A(V,G,\Do)$ to be pointed.

\medbreak

We shall keep the notation of the previous section. Let $S\subseteq G$ be a set of representative classes of $G/\Gamma$, that is $G=\bigcup_{s\in S} s\Gamma$. For any $s\in S$ define
$$A_s=\{ \overline{v\ot g\#1}: v\in \nic(V), g\in s\Gamma\}.$$

If  $\Do=\{d_{ij}\in \ku: 1\leq i, j\leq \theta, i\neq j\}\cup\{\xi_i\in \ku:1\leq i\leq \theta\}$ is a family of scalars then we define $\mathfrak{R}(\Do,\Gamma)$ as the algebra generated by elements $\mathbb{X}_1, \dots, \mathbb{X}_\theta, \gamma\in\Gamma$ subject to relations
$$\gamma \mathbb{X}_i=\chi_i(\gamma)\, \mathbb{X}_i\gamma, \; \;\; \mathbb{X}^{N_i}_i=\xi_i 1, \;\;\; \mathbb{X}_i\mathbb{X}_j-q_{ij}\,\mathbb{X}_j\mathbb{X}_i=d_{ij} 1.$$
The following result seems to be well-known.
\begin{lema} The algebra $\mathfrak{R}(\Do,\Gamma)$ is basic if and only if $d_{ij}=0=\xi_i$ for all $0\leq i, j\leq\theta$.\qed
\end{lema}

The following result generalizes \cite[Thm. 5.2]{EG1}.

\begin{prop} The following hold:
\begin{enumerate}
  \item[1.]  For any $s\in S$ the space $A_s$ is a subcoalgebra,
  \item[2.]  $A(V,G,\Do)=\oplus_{s\in S}\, A_s$,
  \item[3.]  there is an isomorphism of algebras $A_s^*\simeq \mathfrak{R}(\widehat{\Do}-s\cdot \widehat{\Do},\Gamma)$.
\end{enumerate}
\end{prop}

\pf 1. It follows from the definition of the coproduct of $A(V,G,\Do)$ and the fact that the twist $J_\Do\in \nic(V)\ot \nic(V)$.

3. The product in $A_s^*$ is described as follows. If $X, Y\in A_s^*$, $h\in A_s$ then
\begin{align}\label{tw-product}\langle X\ast Y, h\rangle= \langle X, J^{-1} h j^1\rangle\langle Y, J^{-2} hj^2\rangle,
\end{align}
where $J=j^1\ot j^2=J^1_\Do\ot 1\#1 \ot J^2_\Do\ot 1\#1$, $J^{-1}\ot J^{-2}=J^{-1}$.
A base of $A_s$ is given by $\{\overline{x_1^{r_1}\dots x_\theta^{r_\theta}\ot s\gamma\# 1}:  1\leq r_i \leq N_i, 1\leq i \leq\theta, \gamma\in\Gamma\}$. A base for the dual space $A_s^*$ is given by the family $\{X_1^{r_1}\dots X_\theta^{r_\theta}\gamma:  1\leq r_i \leq N_i, 1\leq i \leq\theta, g\in\Gamma\}$ where
$$\langle X_1^{r_1}\dots X_\theta^{r_\theta}g^* ,  \overline{x_1^{s_1}\dots x_\theta^{s_\theta}\ot s\gamma\# 1}\rangle=\begin{cases} \prod_{i=1}^\theta\; (r_i)!_{q_i} \langle g^* , \gamma\rangle\quad \text{ if }  \; r_i=s_i \forall i,\\
0 \quad \text{ otherwise.} \end{cases}$$
Here $ \langle g^* , g_i\rangle=\chi_i(g)$, for all $g_i\in \Gamma$. It is not difficult to verify that using the product \eqref{tw-product} the following equations hold:
$$g^* \ast f^*=(gf)^*, \; g^* \ast X_i = \chi_i(g)\, X_i g^*,\quad X_i \ast X_j= X_iX_j+ ((\chi_i\chi_j(s)-1)a_{ij}) 1,$$
$$ X_i \ast X_i^{N_i-1}=(\xi_i-\chi_i^{N_i}(s)\xi_i) 1, \quad  X_i^{l}\ast X_i^{k}= X_i^{k+l}\,\, \text{for all } \,  k+l<N_i-1,$$
for all $g,f\in\Gamma$, $ 1\leq i,j\leq \theta$. We shall give the proof of third equality, the proofs of the other equations are done in a completely similar way. Let $g\in \Gamma$ then
\begin{align*} \langle X_i \ast X_j, \overline{1\ot sg\# 1}\rangle&=\langle X_i, \overline{J^{-1}_\Do j^1_{s\cdot\Do}\ot sg\# 1}\rangle\langle X_j, \overline{J^{-2}_\Do j^2_{s\cdot\Do}\ot sg\# 1}\rangle\\
&= - a_{ij} \langle X_i, \overline{x_i \ot sg\# 1}\rangle \langle X_j, \overline{x_j \ot sg\# 1}\rangle +\\
 &\qquad +\chi_i\chi_j(s) a_{ij} \langle X_i, \overline{x_i \ot sg\# 1}\rangle \langle X_j, \overline{x_j \ot sg\# 1}\rangle\\
 &= - a_{ij} +\chi_i\chi_j(s) a_{ij}.
\end{align*}
For the second equation we are using that the  coefficient of the term $x_i\ot x_j$ of  $\exp_{q_{ij}}(B_{ij})^{-1}$ is $-a_{ij}$. Since $\langle X_i \ast X_j, \overline{x_ix_j \ot 1 \# 1}\rangle=1$ and
$\langle X_i \ast X_j, \overline{v \ot 1 \# 1}\rangle=0$ for any $v\in \nic(V)$ different from $1$ and $x_ix_j$ then the result follows.

Observe that since $X_iX_j=q_{ij}\, X_jX_i$ then
$$ X_i \ast X_j- q_{ij}\, X_j \ast X_i=(q_{ij} a_{ji}-a_{ij}+ \chi_i\chi_j(s)a_{ji}-q_{ij}\chi_i\chi_j(s)a_{ij} ) 1. $$
Whence there is a well-defined projection $\mathfrak{R}(\widehat{\Do}-s\cdot \widehat{\Do},\Gamma)\to A_s^*$ and since both algebras have the same dimension they must be isomorphic.
\epf

We can generalize \cite[Corollary 5.3]{EG1}.
\begin{cor} The Hopf algebra $A(V,G,\Do)$ is pointed if and only if $\widehat{\Do}$ is $G$-invariant. In particular if $\Do$ is q-symmetric then $A(V,G,\Do)$ is pointed if and only if $\Do$ is $G$-invariant if and only if $A(V,G,\Do)\simeq \nic(V)\#\ku G.$
\end{cor}
\pf $A(V,G,\Do)$ is pointed if and only if $A_s$ is pointed for all $s\in S$ if and only if $A^*_s$ is basic for all $s\in S$ if and only if $\mathfrak{R}(\widehat{\Do}-s\cdot \widehat{\Do},\Gamma,s)$ is basic for all $s\in S$ if and only if $\widehat{\Do}$ is $s$-invariant for all $s\in S$.

If $\Do$ is q-symmetric then $a_{ij}-q_{ij} a_{ji}=2 a_{ij}$ thus $\widehat{\Do}$ is $G$-invariant if and only if  $\Do$ is $G$-invariant. If $\Do$ is $G$-invariant then by Remark \ref{rmk-coprod} $A(V,G,\Do)\simeq \nic(V)\#\ku G$ thus $A(V,G,\Do)$ is pointed.

\epf

\begin{question}  If $\widehat{\Do}$ is $G$-invariant then $A(V,G,\Do)\simeq \nic(V)\#\ku G$ ?
\end{question}

As a last remark I would like to point out that the dual of the Hopf algebra $A(V,G, W, F, J, \Do)$ has coradical $\ku^G$ and this family could be helpful to the study of Hopf algebras with coradical a Hopf subalgebra that has recently began in \cite{AV}.

\end{document}